\documentclass[12pt]{amsart}
\usepackage[latin1]{inputenc}
\usepackage[english]{babel}
\usepackage{amsmath,amssymb,amsthm,amsfonts, color}
\usepackage{latexsym,dsfont}
\usepackage{enumerate}
\usepackage{graphicx}
\usepackage{latexsym}
\usepackage[latin1]{inputenc}
\usepackage[english]{babel}
\usepackage{amsmath,amssymb,amsthm,amsfonts}
\usepackage[all]{xy}
\usepackage{fancyhdr}
\makeatletter
\@namedef{subjclassname@2010}{ \textup{2010} Mathematics Subject Classification}
\makeatother

\newtheorem{theorem}{Theorem}%[section]

\newtheorem{lemma}{Lemma}%[section]

\theoremstyle{definition}
\newtheorem{definition}[lemma]{Definition}

\newtheorem{remark}[lemma]{Remark}

\numberwithin{equation}{section}

\numberwithin{equation}{section}

\numberwithin{equation}{section}
\newcommand{\norm}[1]{\left\lVert#1\right\rVert}

\begin{document}
\title{Pointwise mutipliers of Orlicz function spaces and factorization}
%\thanks{{\rm *}This publication has been produced during scholarship 
%period of the first author at the Lule{\aa} University of Technology, thanks to a Swedish Institute scholarschip 
%(number 0095/2013).}

\author[Le\'snik]{Karol Le\'snik}
\address[Karol Le{\'s}nik]{Institute of Mathematics\\
 Pozna\'n University of Technology, ul. Piotrowo 3a, 60-965 Pozna{\'n}, Poland}
\email{\texttt{klesnik@vp.pl}}
\author[Tomaszewski]{Jakub Tomaszewski}
\address[Jakub Tomaszewski]{Institute of Mathematics\\
Pozna\'n University of Technology, ul. Piotrowo 3a, 60-965 Pozna{\'n}, Poland}
\email{\texttt{uzy93v11@gmail.com}}
\maketitle

\vspace{-9mm}

\begin{abstract}
In the paper we find representation of the space of pointwise multipliers between two Orlicz function spaces, which appears to be another Orlicz space and the formula for the Young function generating this space is given. %Thus we prove the conjecture posted in \cite{KLM14}. 
Further, we apply this result to find necessary and sufficient conditions for factorization of Orlicz function spaces. 

\end{abstract}

\renewcommand{\thefootnote}{\fnsymbol{footnote}}


\section{Introduction}

The space of pointwise multipliers $M(L^{\varphi_1},L^{\varphi})$ is the space of all functions $x$, such that $xy\in L^{\varphi}$ for each $y\in L^{\varphi_1}$, equipped with the operator norm.  
The problem of identifying such spaces was investigated by many authors, starting from Shragin \cite{Sh57}, Ando \cite{An60}, O'Neil \cite{ON65} and Zabreiko-Rutickii \cite{ZR67}, who gave a number of partial answers.   

These investigations were continued in number of directions and results were presented in different forms. One of them is the following result from Maligranda-Nakaii paper \cite{MN10}, which states that if for two given Young functions  $\varphi,\varphi_1$ there is a third one $\varphi_2$ satisfying
\begin{equation}\label{odwr}
\varphi_1^{-1}\varphi_2^{-1}\approx \varphi^{-1},
\end{equation}
then 
\[
M(L^{\varphi_1},L^{\varphi})=L^{\varphi_2}.
\]

This result, however, neither gives any information when such a function $\varphi_2$ exists, nor says anything  how to find it. Further, it was proved in \cite{KLM13} that condition (\ref{odwr}) is  necessary for a wide class of $\varphi,\varphi_1$ functions satisfying some additional properties, but at the same time Example 7.8 from \cite{KLM13} ensures that in general it is not a case, i.e. there are functions $\varphi,\varphi_1$ such that no Young function $\varphi_2$ satisfies (\ref{odwr}), while
\[
M(L^{\varphi_1},L^{\varphi})=L^{\infty}.
\]

On the other hand, there is a natural candidate for function $\varphi_2$ satisfying
\[
M(L^{\varphi_1},L^{\varphi})=L^{\varphi_2}.
\]
Such a function is the following generalization of Young conjugate function (a kind of generalized Legendre transform considered also in convex analysis, for example in \cite{St96}) defined for two Orlicz functions $\varphi,\varphi_1$ as 
\[
\varphi\ominus \varphi_1(t)=\sup_{s>0}\{\varphi(st)-\varphi_1(s)\}.
\]
The function $\varphi\ominus \varphi_1$ is called to be conjugate to $\varphi_1$ with respect to $\varphi$. 
 
Also  in \cite{KLM13}  this construction was compared with condition (\ref{odwr}) and it happens that very often $\varphi_2=\varphi\ominus \varphi_1$ satisfies  (\ref{odwr}), but once again Example 7.8 from \cite{KLM13} shows that, in general,  $\varphi_2=\varphi\ominus \varphi_1$ need not satisfy  (\ref{odwr}). In this example, anyhow, there holds $L^{\infty}=L^{\varphi\ominus \varphi_1}$, so that $M(L^{\varphi_1},L^{\varphi})=L^{\varphi\ominus \varphi_1}$. 
Therefore, it is natural to expect that in general
\begin{equation}\label{con}
M(L^{\varphi_1},L^{\varphi})=L^{\varphi\ominus \varphi_1}.
\end{equation}

In fact, this was already stated for N-functions by Maurey in \cite{Ma74}, but his proof depends heavily on the false conjecture, that the construction $\varphi\ominus \varphi_1$ enjoys involution property, i.e. $\varphi\ominus(\varphi\ominus \varphi_1)=\varphi_1$ (see Example 7.12 in \cite{KLM13} for counterexample). 

On the other hand, the conjecture (\ref{con}) was already proved for Orlicz sequence spaces by Djakov and Ramanujan in \cite{DR00}, where they used a slightly modified construction $\varphi\ominus \varphi_1$ (the supremum is taken only over $0<s\leq 1$). This modification appeared to be appropriate for sequence case, because then only behaviour of Young functions for small arguments is important, while cannot be used for function spaces. Anyhow, we will borrow some ideas from \cite{DR00}. 

In our main Theorem \ref{main} we prove that (\ref{con}) holds in full generality for Orlicz function spaces, as well over finite and infinite measure. 
Then we use this result to find that $\varphi_2=\varphi\ominus \varphi_1$ satisfies (\ref{odwr}) if and anly if $L^{\varphi_1}$ factorizes $L^{\varphi}$, which completes the discussion from \cite{KLM14}.

%%%%%%%%%%%%%%%%%%%%%%%%%%%%%%%%%%%%%%%%%%%%%%
\section{Notation and preliminaries}

Let $L^0=L^0(\Omega,\Sigma,\mu)$ be the space of all classes of $\mu$-measurable, real valuable functions on $\Omega$, where $(\Omega,\Sigma,\mu)$ is a $\sigma$-finite complete measure space. A Banach space $X\subset L^0$ is called the {\it  Banach ideal space} if it satisfies the so called ideal property, i.e. $x\in L^0,y\in X$ with $|x|\leq |y|$ implies $x\in X$ and $\|x\|_X\leq \|y\|_X$ (here $|x|\leq |y|$ means that $|x(t)|\leq |y(t)|$ a.e.), and it contains a {\it weak unity}, i.e. a function $x\in X$ such that $x(t)>0$ for $\mu$-a.e. $t\in \Omega$. When $(\Omega,\Sigma,\mu)$ is purely nonatomic measure spaces, the respective space is called {\it Banach function space} (abbreviation B.f.s.), while in case of $\mathbb{N}$ with counting measure we shall speak about {\it Banach sequence space}.
A Banach ideal space $X$ satisfies the {\it Fatou property} when given a sequence $(x_n)\subset X$, satisfying $x_n\uparrow x$ $\mu$-a.e. and $\sup_n\|x_n\|_X<\infty$, there holds $x\in X$ and $\|x\|_X\leq \sup_n\|x_n\|_X$. 

Writing $X=Y$ for two B.f.s. we mean that they are equal as set, but norms are just equivalent. Recall also that for Banach ideal spaces $X,Y$ the inclusion $X\subset Y$ is always continuous, i.e. there is $c>0$ such that $\|x\|_Y\leq c\|x\|_X$ for each $x\in X$. 
 
For two given Banach ideal spaces $X,Y$ over the same measure space $(\Omega,\Sigma,\mu)$, the space of pointwise multipliers from $X$ to $Y$ is defined as
\[
M(X,Y)=\{y\in L^0: xy\in Y\ {\rm for\ all\ }y\in X\}
\]
with the natural operator norm 
\[
\|y\|_{M(X,Y)}=\sup_{\|x\|_X\leq1}\|xy\|_Y. 
\]
Such a space may be trivial, for example $M(L^p,L^q)=\{0\}$ when $p>q$, and therefore it need not be a Banach ideal space in the sense of above definition. Anyhow, it is a Banach space with the ideal property (see for example \cite{MP89}). When there is no risk of confusion we will just write $\|\cdot\|_M$ for the norm of $M(X,Y)$.

A function $\varphi:[0,\infty)\to [0,\infty]$ will be called a Young function if it is convex, non-decreasing and $\varphi(0)=0$. We will need the following parameters
\[
a_{\varphi}=\sup\{t\geq 0:\varphi(t)=0\} {\rm\ and\ } b_{\varphi}=\sup\{t\geq 0:\varphi(t)<\infty\}.
\]
Let $\varphi$ be a Young function. The {\it Orlicz space} $L^{\varphi}$ is defined as
\[
L^{\varphi}=\{x\in L^0:I_{\varphi}(\lambda x)<\infty {\rm\ for\ some}\ \lambda>0\},
\]
where the modular $I_\varphi$ is given by
\[
I_{\varphi}(x)=\int_{\Omega}\varphi(|x|)d\mu
\]
and the Luxemburg-Nakano norm is defined as
\[
\|x\|_{\varphi}=\inf\{\lambda>0:I_{\varphi}( \frac{x}{\lambda})\leq 1\}.
\]
We point out here that the function $\varphi\equiv 0$ is excluded from the definition of Young functions, but we allow $\varphi(u)=\infty$ for each $u>0$ and understand that in this case $L^{\varphi}=\{0\}$. 

We will often use the following relation between norm and modular. For $x\in L^{\varphi}$ 
\begin{equation}\label{nm}
\|x\|_{\varphi}\leq 1 \Rightarrow I_{\varphi}(x)\leq\norm{x}_{\varphi},
\end{equation}
(see for example \cite{Ma89}).

%We start with some comments on construction $\varphi\ominus \varphi_1$. 

For a given two Young functions $\varphi,\varphi_1$ let us define the mentioned construction of another Young function $\varphi\ominus \varphi_1$, this is 
\[
\varphi\ominus \varphi_1(u)=\sup_{0\leq s}\{\varphi(su)-\varphi_1(s)\}.
\]
Notice that it is a natural generalization of conjugate function in a sense of Young, i.e. $\varphi\ominus \varphi_1$ is called the conjugate function to $\varphi_1$ with respect to $\varphi$. Of course, when $\varphi(u)=u$ we get just the classical conjugate function $\varphi_1^*$ to $\varphi_1$. 
In the above definition, one may be confused by possibility of appearance of indefinite symbol $\infty-\infty$ when $b_{\varphi},b_{\varphi_1}<\infty$. To avoid such a situation we understand that the supremum is taken over $0<s<b_{\varphi_1}$ when $b_{\varphi_1}<\infty$ and $\varphi_1(b_{\varphi_1})=\infty$, or over $0<s\leq b_{\varphi_1}$ when $b_{\varphi_1}<\infty$, but $\varphi_1(b_{\varphi_1})<\infty$. 
Of course, functions $\varphi,\varphi_1$  and $\varphi\ominus \varphi_1$ satisfy the generalized Young inequality, i.e. 
\[
\varphi(uv)\leq \varphi\ominus \varphi_1(u)+\varphi_1(v)
\]
for each $u,v\geq 0$.

We will  also need the following construction.
\begin{definition} For two Young functions $\varphi,\varphi_1$ and $0<a\leq b_{\varphi_1}$ we define
\[
\varphi\ominus_a \varphi_1(u)=\sup_{0\leq s\leq a}\{\varphi(su)-\varphi_1(s)\}.
\]
\end{definition}

Such defined function $\varphi\ominus_a \varphi_1$ enjoys the following elementary properties.

\begin{lemma} Let $\varphi,\varphi_1$ be two Young functions. 
\begin{enumerate}

\item $\varphi\ominus_a \varphi_1$ is Young function for each $0<a\leq b_{\varphi_1}$.

%\item If $b_{\varphi}, b_{\psi}<\infty$ then $b_{\varphi\ominus\psi}=\frac{b_{\varphi}}{b_{\psi}}.$

\item For each $t\geq 0$ there holds
$$
\lim\limits_{a\rightarrow b_{\varphi_1}^{-}}\varphi\ominus_a \varphi_1(u)=\varphi\ominus \varphi_1(u).
$$
 \end{enumerate}
\end{lemma}

\begin{remark}\label{rem1}
Notice that dilations of Young functions do not change Orlicz spaces, i.e. when $\varphi$ is a Young function and $\psi$ is defined by $\psi(u)=\varphi(au)$ for some $a>0$, then $L^{\varphi}=L^{\psi}$. It gives a reason to expect that dilating  $\varphi,\varphi_1$ results in dilation of $\varphi\ominus\varphi_1$. In fact, let $\varphi,\varphi_1$ be Young functions and put $\psi(u)=\varphi(au)$, $\psi_1(u)=\varphi_1(bu)$. Then 
\[
\psi\ominus\psi_1(u)=\sup_{0<s}(\varphi(aus)-\varphi_1(bs))=\sup_{0<s}(\varphi(aus/b)-\varphi_1(s))=\varphi\ominus\varphi_1(au/b).
\]
Moreover, if  $b_{\varphi}=b_{\varphi_1}<\infty$, then supremum in the definition of $\varphi\ominus\varphi_1$ is attained for each $u<1$, i.e. for each $u<1$ there is $0<s<b_{\varphi_1}$ such that  $\varphi\ominus\varphi_1(u)=\varphi(us)-\varphi_1(s)$. In particular, $b_{\varphi\ominus\varphi_1}=1$. 
\end{remark}

\begin{remark} \label{rem2}
Let us also recall that a fundamental function $f_{\varphi}$ of an Orlicz space $L^{\varphi}$ is given by the formula $f_{\varphi}(t)=\frac{1}{\varphi^{-1}(1/t)}$,
for $0<t<\mu(\Omega)$ and $f_{\varphi}(0)=0$, where by $\varphi^{-1}$ we understand the right continuous inverse of $\varphi$, i.e. 
$\varphi ^{-1}(v) = \inf \{u\geq 0: \varphi (u)>v\}$ (more informations about fundamental functions of symmetric spaces may be found in \cite{BS88}) .  In particular, the fundamental function of $L^{\varphi}$ is right-continuous at $0$  if and only if $b_{\varphi}=\infty$, or equivalently, $b_{\varphi}=\infty$ if and only if for each $\varepsilon>0$ there is $\delta>0$ such that if $A\in \Sigma$, $\mu(A)<\delta$ then $\|\chi_{A}\|_{\varphi}<\varepsilon$. %In another words, the fundamental function of $L^{\varphi}$ is right-continuous at $0$. Moreover, if  $b_{\varphi}<\infty$, then 
\end{remark}

%%%%%%%%%%%%%%%%%
\section{Multipliers of Orlicz function spaces}

\begin{lemma}\label{lemzero}
Let $\varphi,\varphi_1$ be Young functions such that $b_{\varphi}<\infty$ and $b_{\varphi_1}=\infty.$ Then 
$$
M(L^{\varphi_1},L^{\varphi})=\{0\}.
$$
\end{lemma}
\proof The proof follows immediately from Proposition 3.2 in \cite{KLM13}, since under our assumptions $L^{\varphi_1}\not \subset L^\infty$ but  $L^{\varphi}\subset L^\infty$. 
\endproof

\begin{lemma}\label{factobaskzanurz}
Let $\varphi,\varphi_1$ be Young functions and $b_{\varphi}<\infty$. Then 
$$
M(L^{\varphi_1},L^{\varphi})\subset L^{\infty}.
$$  
\end{lemma}
\proof
Suppose that $M(L^{\varphi_1},L^{\varphi})\not\subset L^{\infty}$. Then there exists $0\leq y\in M(L^{\varphi_1},L^{\varphi})$ such that $\|y\|_{M}=1$ and for each $n>0$
$$
\mu(\{t\in \Omega: y(t)\geq n\})>0.
$$ 
Denote $A_n=\{t\in \Omega: y(t)\geq n\}$ for $n\in\mathbb{N}$. Then $\norm{n\chi_{A_n}}_{M}\leq 1$ and for $A_{n_0}$ chosen such that $\mu(A_{n_0})<\infty$, it follows for $n>n_0$
$$
\|y\|_{M}\geq\|n\chi_{A_n}\|_{M}\geq \frac{n}{\|\chi_{A_{n_0}}\|_{\varphi_1}}\|\chi_{A_n}\chi_{A_{n_0}}\|_{\varphi}=\frac{n}{\|\chi_{A_{n_0}}\|_{\varphi_1}}\|\chi_{A_n}\|_{\varphi}\geq \frac{nb_{\varphi}^{-1}}{\|\chi_{A_{n_0}}\|_{\varphi_1}}.
$$ 
This contradiction shows that $M(L^{\varphi_1},L^{\varphi})\subset L^{\infty}$.
\endproof

We are in a position to prove the main theorem. 

\begin{theorem}\label{main}
Let $\varphi,\varphi_1$ be Young functions. Then
$$
M(L^{\varphi_1},L^{\varphi})=L^{\varphi\ominus\varphi_1}.
$$
\end{theorem}
\proof
The inclusion
\begin{equation}
L^{\varphi\ominus\varphi_1}\subset M(L^{\varphi_1},L^{\varphi}) \label{embedding}
\end{equation}
 is well known (see  \cite{An60}, \cite{KLM13}, \cite{MN10} or \cite{ON65}) and follows from equivalence of generalized Young inequality and inequality $\varphi_1^{-1}(\varphi\ominus\varphi_1)^{-1}\lesssim \varphi^{-1}$. For the completeness of presentation we present the proof which employs the generalized Young inequality  directly. 
If $\varphi\ominus\varphi_1(u)=\infty$ for each $u>0$ then $L^{\varphi\ominus\varphi_1}=\{0\}$ and inclusion trivially holds. Suppose $L^{\varphi\ominus\varphi_1}\neq\{0\}$, i.e. $\varphi\ominus\varphi_1(u)<\infty$ for some $u>0$.
% Note that by definition of $\varphi\ominus\varphi_1$, the generalized Young inequality
%$$
%\varphi\left(st\right)\leq\varphi\ominus\varphi_1\left(s\right)+\varphi_1\left(t\right)
%$$
%is satisfied for all $s,t>0$ %such that $\varphi(st)<\infty$
% (sprawdic co z bfi).
Let $y\in L^{\varphi\ominus\varphi_1}$ and $x\in L^{\varphi_1}$ be such that
$$
\norm{y}_{\varphi\ominus\varphi_1}\leq\frac{1}{2}\\ \text{ and }\\ \norm{x}_{\varphi_1}\leq\frac{1}{2}.
$$
Then generalized Young inequality gives 
$$
 I_{\varphi}(y x)\leq I_{\varphi\ominus\varphi_1}(y)+I_{\varphi_1}(x)\leq 1.
$$ 
Consequently $y x\in L^{\varphi}$ and $\norm{y x}_{\varphi}\leq 1$. Therefore, $L^{\varphi\ominus\varphi_1}\subset M(L^{\varphi_1},L^{\varphi})$ and
$$
\norm{y}_{M}\leq 4\norm{y}_{\varphi\ominus\varphi_1}.
$$

To prove the second inclusion it is enough to indicate a constant $c>0$  such that for each simple function $y\in M(L^{\varphi_1},L^{\varphi})$ there holds
\begin{equation}\label{fp}
\|y\|_{\varphi\ominus\varphi_1}\leq c\|y\|_M.
\end{equation}
In fact, it follows directly from the Fatou property of both $L^{\varphi\ominus\varphi_1}$ and $M(L^{\varphi_1},L^{\varphi})$ spaces (it is elementary fact that $M(X,Y)$ has the Fatou property when $Y$ has so). Let $0\leq y\in M(L^{\varphi_1},L^{\varphi})$ and $0\leq y_n\uparrow y$ $\mu$-a.e., where $y_n$ are simple functions. Then, by (\ref{fp}),
\[
\|y_n\|_{\varphi\ominus\varphi_1}\leq c\|y_n\|_M\rightarrow c\|y\|_M
\]
and so the  Fatou property of $L^{\varphi\ominus\varphi_1}$ implies $y\in L^{\varphi\ominus\varphi_1}$ and  $\|y\|_{\varphi\ominus\varphi_1}\leq c\|y\|_M$.

The proof of (\ref{fp}) will be divided into four cases, depending on finiteness of $b_{\varphi}$ and $b_{\varphi_1}$.

%In rest of the proof we can assume, without lost of generality, that $$M(L^{\varphi_1},L^{\varphi})\neq\{0\}.$$ 

Consider firstly the most important case $b_{\varphi}=b_{\varphi_1}=\infty$.
Let $0\leq y\in M(L^{\varphi_1},L^{\varphi})$ be a simple function of the form 
$y=\sum_k a_k\chi_{B_k}$ and such that  $\norm{y}_{M}\leq\frac{1}{2}$. We will show that for each $a>1$ 
$$
I_{\varphi\ominus_a\varphi_1}(y)\leq 1.
$$ 
Let $a>1$ be arbitrary. For each $a_k$ there exists $b_k\geq 0$ such that
$$
\varphi(a_kb_k)=\varphi\ominus_a\varphi_1(a_k)+\varphi_1(b_k).
$$ 
This is, for $x=\sum_k b_k\chi_{B_k}$, there holds  $\varphi(xy)=\varphi\ominus_a\varphi_1(x)+\varphi_1(y)$. Note that from definition of $\varphi\ominus_a\varphi_1$ we have  $x(t)\leq a$ for each $t\in \Omega$. Further, since $b_{\varphi_1}=\infty$, there exists $t_a>0$ such that $\|\chi_{A}\|_{\varphi_1}\leq\frac{1}{a}$ for each $A\subset \Omega$ with $\mu(A)<t_a$ (see Remark \ref{rem2}). Suppose $\mu(\Omega)=\infty$. 
Since $(\Omega,\Sigma,\mu)$ is $\sigma$-finite and atomless, we can divide $\Omega$ into a sequence of pairwise disjoint sets $(A_n)$ with $\mu(A_n)= t_a$ for each $n\in \mathbb{N}$ and $\Omega =\bigcup A_n$. In the case of $\mu(\Omega)<\infty$ the sequence $(A_n)$ may be chosen finite and such that $\mu(A_n)=\delta\leq t_a$ for each $n=1,\dots,N$ with $\Omega =\bigcup A_n$.
 
In any case, for $A_n$ we have 
$$
\|y x\chi_{A_n}\|_{\varphi}\leq \|y\|_M\|x\chi_{A_n}\|_{\varphi_1}\leq\frac{a}{2}\|\chi_{A_n}\|_{\varphi_1}\leq\frac{1}{2},
$$ 
because $\mu(A_n)\leq t_a$ and $x(t)\leq a$ for $t\in \Omega$. In consequence, using inequality  $\varphi_1(x)\leq\varphi(yx)$,  we have for each $A_n$
\begin{equation}\label{equ1}
I_{\varphi_1}(x\chi_{A_n})\leq I_{\varphi}(y x\chi_{A_n})\leq\|y x\chi_{A_n}\|_{\varphi}\leq\frac{1}{2}.
\end{equation}
%indukcja
Define now 
$$
x_n=\sum_{k=1}^nx\chi_{A_k}.
$$ 
We claim that $I_{\varphi_1}(x_n)\leq\frac{1}{2}$ for each $n$. It will be shown by induction. 
 For $n=1$ it comes from (\ref{equ1}). Let $n>1$ and suppose 
$$
I_{\varphi_1}(x_{n-1})\leq\frac{1}{2}.
$$ 
It follows
 $$
I_{\varphi_1}(x_{n})=I_{\varphi_1}(x_{n-1})+I_{\varphi_1}(x\chi_{A_n})\leq 1,
$$ 
hence $\norm{x_n}_{\varphi_1}\leq 1$. Moreover,
$$
\norm{y x_n}_{\varphi}\leq\frac{1}{2}\norm{x_n}_{\varphi_1}\leq\frac{1}{2}
$$
together with inequality $\varphi_1(x)\leq\varphi(yx)$ imply  
$$
I_{\varphi_1}(x_{n})\leq I_{\varphi}(y x_{n})\leq\norm{y x_n}_{\varphi}\leq\frac{1}{2}.
$$ 
It means we proved the claim and can proceed with the proof. %In case $I=[0,1]$ there exist $n_0\in\mathbb{N}$ such $$x=x_{n_0},$$ therefore $$\norm{x}_{\psi}=\norm{x_{n_0}}_{\psi}\leq 1.$$
  %In case $I=[0,\infty)$ 
Clearly, $x_n\uparrow x$ $\mu$-a.e., thus from the Fatou property of $L^{\varphi_1}$ we obtain  that $x\in L^{\varphi_1}$ and 
$$
\|x\|_{\varphi_1}\leq\sup\limits_n\norm{x_n}_{\varphi_1}\leq 1.
$$  
Finally, inequality $\varphi\ominus_a\varphi_1(y)\leq\varphi(yx)$
together with $\|y x\|_{\varphi}\leq\frac{1}{2}\|x\|_{\varphi_1}\leq\frac{1}{2}$ give
 $$
I_{\varphi\ominus_a\varphi_1}(y)\leq I_{\varphi}(y x)\leq \norm{y x}_{\varphi}\leq\frac{1}{2}.
$$
Applying Fatou Lemma we obtain 
 $$
I_{\varphi\ominus\varphi_1}(y)= \int\varphi\ominus\varphi_1(y)d\mu\leq\liminf_{a\rightarrow\infty}\int\varphi\ominus_a\varphi_1(y)d\mu\leq\frac{1}{2}.
$$ 
In consequence $y\in L^{\varphi\ominus\varphi_1}$ with $\|y\|_{\varphi\ominus\varphi_1}\leq 1$. This gives also constant for inclusion, i.e. 
$$
\|y\|_{\varphi\ominus\varphi_1}\leq 2\|y\|_{M},
$$ 
when $y\in M(L^{\varphi_1},L^{\varphi})$.

 Let us consider the second case, this is $b_{\varphi}=\infty$ and $b_{\varphi_1}<\infty$. Without loss of generality we can assume that $b_{\varphi_1}>1$ (see Remark \ref{rem1}). Let $0\leq y\in M(L^{\varphi_1},L^{\varphi})$ 
be a simple function satisfying $\|y\|_{M}\leq\frac{1}{2b_{\varphi_1}}$. Notice that  $b_{\varphi}=\infty$ with $b_{\varphi_1}<\infty$ imply that $b_{\varphi\ominus\varphi_1}=\infty$. Moreover, as before, there exists a simple function $x$ such that $0<x(t)\leq b_{\varphi_1}$ for each $t\in \Omega$ and 
$$
\varphi(yx)=\varphi\ominus\varphi_1(y)+\varphi_1(x)
$$ 
(see Remark \ref{rem1}). 
As before, we can find $t_0>0$ such that $\mu(A)<t_0$ implies $\|\chi_{A}\|_{\varphi_1}\leq 1$. Selecting  the sequence $(A_n)$ like previously, but with  $\mu(A_n)\leq t_0$ for each $A_n$, we obtain 
$$
\norm{y x\chi_{A_n}}_{\varphi}\leq \frac{b_{\varphi_1}}{2b_{\varphi_1}}\|\chi_{A_n}\|_{\varphi_1}\leq\frac{1}{2}.
$$ 
%is chosen to satisfy $\mu(A_n)= t_0$ for each $n\in \mathbb{N}$ and $\Omega =\bigcup A_n$, when  $\mu(\Omega)=\infty$, or $\mu(A_n)=\delta\leq t_0$ for each $n=1,\dots,N$, $\Omega =\bigcup A_n$, when $\mu(\Omega)<\infty$. 
Define further
$$
x_n=\sum_{k=1}^nx\chi_{A_k}.
$$ 
Then it may be proved by the same induction as before, that  $I_{\varphi_1}(x_n)\leq\frac{1}{2}$ for each $n$. Following respective steps from previous case we get
$$
\|y\|_{\varphi\ominus\varphi_1}\leq 2b_{\varphi_1}\norm{y}_{M}.
$$
 
Let now $b_{\varphi},b_{\varphi_1}<\infty$.
We can assume that $b_{\varphi_1}=b_{\varphi}=1$ (see Remark \ref{rem1}). From Lemma \ref{factobaskzanurz} it follows that there exists a constant $c\geq 1$ such that for each $y\in M(L^{\varphi_1},L^{\varphi})$ we have 
$$
\|y\|_{\infty}\leq c\|y\|_{M}.
$$
 Let $0\leq y\in M(L^{\varphi_1},L^{\varphi})$ be a simple function and $\|y\|_{M}\leq\frac{1}{4c}$.
We have  $y(t)\leq \frac{1}{4c}\leq b_{\varphi\ominus\varphi_1}$ (cf. Remark \ref{rem1}) for almost every $t\in \Omega$, therefore $\varphi\ominus\varphi_1(y(t))<\infty$.
Consequently, we can choose a simple function $x$ satisfying   
$$
\varphi(yx)=\varphi\ominus\varphi_1(y)+\varphi_1(x).
$$ 
Then $x(t)\leq b_{\varphi}=1$ for each $t\in \Omega$. Further, we can find $t_0>0$ so that inequality 
$$
\|\chi_{A}\|_{\varphi_1}\leq 2
$$
 is fulfilled for each $A$ with $\mu(A)\leq t_0$, just because $\lim_{t\to 0^+} f_{\varphi}(t)=b_{\varphi}=1$.  Choosing a sequence $(A_n)$ as in previous cases we get
$$
\|y x\chi_{A_n}\|_{\varphi}\leq \frac{1}{4c}\|\chi_{A_n}\|_{\varphi_1}\leq\frac{1}{2}.
$$
Once again we can show by induction that for each $x_n=\sum_{k=1}^nx\chi_{A_k}$ there holds  $I_{\varphi_1}(x_n)\leq\frac{1}{2}$. Therefore $\|x_n\|_{\varphi_1}\leq 1$ and, by the Fatou property of $L^{\varphi_1}$,  $\|x\|_{\varphi_1}\leq 1$. It follows
$$
\|y x\|_{\varphi}\leq1
$$ 
and by inequality $\varphi\ominus\varphi_1(y)\leq\varphi(yx)$ we obtain
$$
I_{\varphi\ominus\varphi_1}(y)\leq I_{\varphi}(y x)\leq\|y x\|_{\varphi}\leq 1.
$$
Consequently
$$
\|y\|_{\varphi\ominus\varphi_1}\leq 4c\|y\|_{M}.
$$

Finally, there left the trivial case of $b_{\varphi}<\infty$, $b_{\varphi_1}=\infty$ to consider. However, Lemma \ref{lemzero} with the embedding (\ref{embedding}) give
$$
L^{\varphi\ominus\varphi_1}=M(L^{\varphi_1},L^{\varphi})=\{0\}
$$
and the proof is finished.
\endproof

\section{Factorization}

Recall that given two B.f.s. $X,Y$ over the same measure space, we say that $X$ factorizes $Y$ when 
\[
X\odot M(X,Y)=Y,
\]
where 
\[
X\odot M(X,Y)=\{z\in L^0:z=xy {\rm\ for\ some\ }x\in X, y\in M(X,Y)\}.
\]
The idea of such factorization goes back to Lozanovskii, who proved that each B.f.s. factorizes $L^1$. For more informations on factorization and its importance we send a reader to papers \cite{CDS08}, \cite{KLM14} and \cite{Sc10}  which are devoted mainly to this subject. 

Also in \cite{KLM14} one may find a discussion on factorization of Orlicz spaces (and even more general Calder\'on-Lozanovskii spaces). Having in hand our representation $M(L^{\varphi_1},L^{\varphi})= L^{\varphi\ominus\varphi_1}$ we are able to complete this discussion by proving sufficient and necessary conditions for factorization in terms of respective Young functions. 

%that by $\varphi^{-1}$ we understand the right continuous inverse of $\varphi$, i.e. 
%$$
%\varphi ^{-1}(v) = \inf \{u\geq 0: \varphi (u)>v\}.
%$$  
We say that equivalence $\varphi_1^{-1}\varphi_2^{-1}\approx \varphi^{-1}$ holds for all [large] arguments when there are constants $c,C>0$ such that 
\[
c\varphi^{-1}(u)\leq \varphi_1^{-1}(u)\varphi_2^{-1}(u)\leq C\varphi^{-1}(u)
\]
for all $u\geq 0$ [for some $u_0> 0$ and all $u>u_0$]. 

\begin{theorem}\label{factorization}
Let $\varphi,\varphi_1$ be two Young functions. Then $L^{\varphi_1}$ factorizes $L^{\varphi}$, i.e.
\[
L^{\varphi_1}\odot M(L^{\varphi_1},L^{\varphi})=L^{\varphi}
\]
 if and only if 
\begin{itemize}
\item[i)] equivalence $\varphi_1^{-1}(\varphi\ominus\varphi_1)^{-1}\approx \varphi^{-1}$ is satisfied for all arguments when $\mu(\Omega)=\infty$.
\item[ii)] equivalence $\varphi_1^{-1}(\varphi\ominus\varphi_1)^{-1}\approx \varphi^{-1}$ is satisfied for large arguments when $\mu(\Omega)<\infty$,
\end{itemize}
\end{theorem}
\proof In the light of Theorem \ref{main} 
\[
L^{\varphi_1}\odot M(L^{\varphi_1},L^{\varphi})=L^{\varphi_1}\odot L^{\varphi\ominus\varphi_1}.
\]
Therefore  $L^{\varphi_1}$ factorizes $L^{\varphi}$ if and only if $L^{\varphi_1}\odot L^{\varphi\ominus\varphi_1}=L^{\varphi}$. The latter, however, is equivalent with $\varphi_1^{-1}(\varphi\ominus\varphi_1)^{-1}\approx \varphi^{-1}$ for all, or for large arguments, depending on $\Omega$, as proved in Corollary 6 from \cite{KLM14}. 
\endproof

%{\color{red} Guestion: is it true that 
%\[
% M( M(L^{\varphi_1},L^{\varphi}),L^{\varphi})=L^{\varphi_1}
%\]
%implies  $\varphi_1^{-1}\varphi_2^{-1}\approx \varphi^{-1}$?
%The inverse implication is known. }


\begin{thebibliography}{99}

\bibitem{An60} T. Ando, {\it On products of Orlicz spaces}, Math.
Ann. {\bf 140} (1960), 174--186. 


%\bibitem{AZ87} J. Appell and P. P. Zabrejko, {\it On the degeneration 
%of the class of differentiable superposition operators in function spaces}, 
%Analysis {\bf 7} (1987),  no. 3-4, 305--312.


%\bibitem{AZ90} J. Appell and P. P. Zabrejko, {\it Nonlinear
%Superposition Operators}, Cambridge University Press, Cambridge 1990. 


\bibitem{BS88} C. Bennett and R. Sharpley, {\it Interpolation of
Operators}, Academic Press, Boston 1988. 


\bibitem{CDS08} J. M. Calabuig, O. Delgado and E. A. S\'anchez P\'erez, 
{\it Generalized perfect spaces}, Indag. Math. (N.S.) {\bf 19} (2008), no. 3, 359--378. 

%\bibitem{SDS08} J. M. Calabuig, O. Delgado, and E. A. S\'anchez P\'erez, {\it Generalized perfect spaces}, Indag. Math. (N.S.) {\bf 19} (3) (2008), 359--378.

\bibitem{DR00} P. B. Djakov and M. S. Ramanujan, {\it Multipliers
between Orlicz sequence spaces}, Turk. J. Math. {\bf 24} (2000),
313--319. 


\bibitem{KLM13}P. Kolwicz, K. Le\'snik and L. Maligranda, {\it
Pointwise multipliers of Calde\-r\'{o}n-Lozanovski\u{\i} spaces, } Math. Nachr.  {\bf 286} (2013), 876--907.

\bibitem{KLM14}{P. Kolwicz, K. Le\'snik and L. Maligranda, \textit{
Pointwise products of some Banach function spaces and factorization }, J. Funct. Anal.  {\bf 266} (2014), no. 2, 616--659.}

%\bibitem{Lo65} G. Ja. Lozanovski{\u \i}, {\it On reflexive spaces generalizing 
%the reflexive space of Orlicz}, Dokl. Akad. Nauk SSSR {\bf 163} (1965), 
%573--576 (in Russian); English transl. in: Soviet Math. Dokl. {\bf 6} (1965), 
%968--971. 

%\bibitem{Lo69} G. Ja. Lozanovski{\u \i}, {\it On some Banach lattices}, Sibirsk. 
%Mat. Zh. {\bf 10} (1969), 584--599 (Russian); English transl. in Siberian Math. 
%J. {\bf 10} (1969), no. 3, 419--431. 


%\bibitem{Lo73} G. Ja. Lozanovski{\u \i}, {\it Certain Banach
%lattices. IV}, Sibirsk. Mat. Zh. {\bf 14} (1973), 140--155 (in Russian); English
%transl. in: Siberian. Math. J. {\bf 14} (1973), 97--108. 


\bibitem{Ma89} L. Maligranda, {\it Orlicz Spaces and Interpolation}, Seminars 
in Mathematics 5, University of Campinas, Campinas SP, Brazil 1989. 


\bibitem{MN10} L. Maligranda and E. Nakai, {\it Pointwise
multipliers of Orlicz spaces}, Arch. Math. {\bf 95} (2010), no. 3, 251--256. 


\bibitem{MP89} L. Maligranda and L. E. Persson, {\it Generalized
duality of some Banach function spaces}, Indag. Math. {\bf 51} (1989),
no. 3, 323--338.


\bibitem{Ma74} B. Maurey, {\it Th\'eor\`emes de factorisation pour les 
op\'erateurs lin\'eaires {\`a} valeurs dans les espaces $L^{p}$}, Ast\'erisque 
{\bf 11} (1974), 1--163.


%\bibitem{Na95} E. Nakai, {\it Pointwise multipliers}, Memoirs of
%the Akashi College of Technology {\bf 37} (1995), 85--94. 
 

\bibitem{ON65} R. O'Neil, {\it Fractional integration in Orlicz
spaces. I}, Trans. Amer. Math. Soc. {\bf 115} (1965), 300--328. 


\bibitem{Re81} S. Reisner, {\it A factorization theorem in Banach 
lattices and its applications to Lorentz spaces}, Ann. Inst. Fourier (Grenoble) 
{\bf 31} (1981), no. 1, 239--255. 



\bibitem{Sc10} A. R. Schep, {\it Products and factors of Banach
function spaces}, Positivity {\bf 14} (2010), 301--319.
 

\bibitem{Sh57} I. V. Shragin, {\it On certain operators in generalized Orlicz spaces}, Dokl. 
Akad. Nauk SSSR (N.S:) {\bf 117} (1957), 40--43 (in Russian).



\bibitem{St96} T. Str\"omberg, {\it The operation of infimal convolution}, Dissertationes Math. (Rozprawy Mat.) {\bf 352} (1996), 58 pp.

\bibitem{ZR67} P. P. Zabreiko and Ja. B. Rutickii, {\it Several remarks on monotone functions} (in Russian), Uch. Zap. Kazan. Gos. Univ. {\bf 127} (1) (1967), 114--126.

\end{thebibliography}
\end{document}